\documentclass[11pt, a4paper]{article}
\usepackage[latin1]{inputenc}
\usepackage[english]{babel}
\usepackage{indentfirst}
\usepackage{amstext}
\usepackage{setspace}
\usepackage{amsfonts}
\usepackage{textcomp}
\usepackage{amssymb}
\usepackage{amscd}
\usepackage{epsf}
\usepackage{graphicx}
\usepackage{epsfig}
\usepackage[mathscr]{eucal}
\usepackage{xr}
\usepackage{verbatim}
\usepackage{color}
\usepackage{shadethm}

\RequirePackage[reqno]{amsmath}

\newcommand {\demo}{\hskip -0.6cm{\bf Proof:  }}

\newcommand{\cqd}{\hfill{$\blacksquare$}}

\newcommand{\A}{\mathrm A}
\newcommand{\E}{\mathrm E}
\newcommand{\X}{\mathrm X}
\newcommand {\N}{\mathbb{N}}
\newcommand {\Z}{\mathbb{Z}}

\newcommand {\G}{\mathrm{G}}
\newcommand {\C}{\mathbb{C}}

\newcommand{\Fo}{\ensuremath{\mathcal{F}_0(\X)}}

\newtheorem{theorem}{Theorem}[section]
\newtheorem{lema}[theorem]{Lemma}
\newtheorem{corolario}[theorem]{Corollary}
\newtheorem{definition}[theorem]{Definition}
\newtheorem{proposition}[theorem]{Proposition}
\newtheorem{example}[theorem]{Example}
\newtheorem{obs}[theorem]{Remark}

\begin{document}

\title{Simplicity of partial skew group rings of abelian groups}
\maketitle
\begin{center}
{\large Daniel Gonçalves}\\
\end{center}  
\vspace{8mm}

\abstract 
Let $\A$ be a ring with local units, $\E$ a set of local units for $\A$, $\G$ an abelian group and $\alpha$ a partial action of $\G$ by ideals of $\A$ that contain local units and such that the partial skew group ring $\A\star_{\alpha} \G$ is associative. We show that $\A\star_{\alpha} \G$ is simple if and only if $\A$ is $\G$-simple and the center of the corner $e\delta_0 (\A\star_{\alpha} \G) e \delta_0$ is a field for all $e\in \E$. We apply the result to characterize simplicity of partial skew group rings in two cases, namely for partial skew group rings arising from partial actions by clopen subsets of a compact set and partial actions on the set level. 
\doublespace


\section{Introduction}

Partial skew group rings are algebraic analogues of C*-partial crossed products and also arise as natural generalizations of skew group rings to the partial action context (see \cite{Ex}, where partial skew group rings are introduced and their associativity studied). 

As it is the case with its C* counterpart, it is important to realize algebras as partial skew group rings (see for example \cite{gonroy}, where Leavitt path algebras are realized as a partial skew group ring and \cite{Boava} where C*-algebras associated to integral domains are realized as a partial crossed product). The idea behind realizing algebras as partial skew group rings is that one can then benefit from the established general theory of partial skew group rings. Nevertheless, general results about partial skew group rings are still underdeveloped, if compared to the abundance of results in the skew group rings or C* partial crossed product context. For example, much of the ideal structure of skew group rings has been described in \cite{Crow, Fisher, Montgomery,Passman} but, to the author´s knowledge, \cite{Ferrero} is the only reference in the literature regarding the ideal structure of partial skew group rings.
In this context, recently \"{O}inert, see \cite{Oinert}, characterized simplicity of skew group rings of Abelian groups. In this paper we generalize  the results in \cite{Oinert} in two ways, namely, to rings with local units and to partial skew group rings. 

Before we proceed we recall some key definitions below. A partial action of a group $\G$ on a set $\Omega$ is a pair $\alpha= (\{D_{t}\}_{t\in \G}, \ \{\alpha_{t}\}_{t\in \G})$, where for each $t\in \G$, $D_{t}$ is a subset of $\Omega$ and $\alpha_{t}:D_{t^{-1}} \rightarrow \Delta_{t}$ is a bijection such that $D_{e} = \Omega$, $\alpha_{e}$ is the identity in $\Omega$, $\alpha_{t}(D_{t^{-1}} \cap D_{s})=D_{t} \cap D_{ts}$ and $\alpha_{t}(\alpha_{s}(x))=\alpha_{ts}(x),$ for all $x \in D_{s^{-1}} \cap D_{s^{-1} t^{-1}}.$ In case $\Omega$ is an algebra or a ring then the subsets $D_t$ should also be ideals and the maps $\alpha_t$ should be isomorphisms. In the topological setting each $D_t$ should be an open set and each $\alpha_t$ a homeomorphism and in the C*-algebra setting each $D_t$ should be a closed ideal and each $\alpha_t$ should be a *-isomorphism.

Associated to a partial action of a group $G$ in a ring $A$ we have the partial skew group ring $\A\star_{\alpha} \G$, which is defined as the set of all finite formal sums $\sum_{t \in G} a_t\delta_t$, where, for all $t \in G$, $a_t \in D_t$ and $\delta_t$ are symbols. Addition is defined in the usual way and multiplication is determined by $(a_t\delta_t)(b_s\delta_s) = \alpha_t(\alpha_{-t}(a_t)b_s)\delta_{t+s}$. For $a=\sum_{t \in G} a_t \delta_t \in A \star _{\alpha} G $, the support of $a$, which we denote by $supp(a)$, is the finite set 
$\{t \in G \, \, : \, \, a_t \neq 0 \}$ and the projection into the $g$ coordinate map, $P_g:  A \star _{\alpha} G  \rightarrow  A$, is given by $P_g \left(\sum_{t \in G} a_t \delta_t \right)=  a_g.$

\section{Simplicity of $\A\star_{\alpha} \G$}

As we mentioned in the introduction, we are particularly interested in rings with local units, not necessarily unital. So from now on we assume that $A$ is a ring with local units, that is, $A$ is a ring such that for every finite set $\{r_1, r_2, ..., r_t\} \subseteq A$ we can find $e \in A$ such that  $e^2=e$ and  $ e r_i=r_i=r_i e$ \ for every  $i \in \{ 1, \cdots , t \}$. Notice that if $E \subseteq A$ is a set of local units for $A$ then $E \delta_0 = \{ e \delta_0: e\in E\}$ is a set of local units for $\A\star_{\alpha} \G$. 


The condition for simplicity of partial skew group rings relies on the definition of $G$-invariant ideals.  This was defined in \cite{Oinert} for skew group rings and in \cite{ELQ} for C* partial crossed products. Below we give the definition adapted to our context, followed by the first lemma in the paper.

\begin{definition} Let $\alpha= (\{D_{t}\}_{t\in \G}, \ \{\alpha_{t}\}_{t\in \G})$ be a partial action of a group $G$ on a ring $A$.  We say that an ideal $I\subseteq A$ is $G-$invariant if $\alpha_g(I \cap D_{-g}) \subseteq I \cap D_g$, for all $g \in G$. If $\{0\}$ and $A$ are the only $G$-invariant ideals of $A$, then we say that $A$ is $G$-simple.
\end{definition}

\begin{lema}\label{lema 1} 
Let $E$ be a set of local units for $A$ and $\alpha=(\{D_t\}_{t\in G}, \{\alpha_t\}_{t \in G} )$ a partial action of an abelian group $G$ such that each ideal $D_t$ has local units. Suppose that $R= A \star _{\alpha} G$ is associative and $A$ is $G$-simple. Then, for each non-zero $r \in A \star _{\alpha}G$, and for each local unit $e\in E$, there exists $r' \in R$ such that:
\begin{enumerate}
\renewcommand{\labelenumi}{\alph{enumi})}
\item $r' \in RrR$
\item $P_0(r')= e$
\item $\# supp(r') \leq \# supp(r) $
\end{enumerate}
\end{lema}
\demo 

Let $r= \sum_{g \in G} a_g \delta_g $ be a non-zero element in $R$. Let $h\in G$ be such that $a_h \neq 0$ and $e_h \in D_h$ be a unit for $a_h$. Notice that $a_h\delta_h \cdot  \alpha_{-h}(e_h) \delta_{-h} = \alpha_h \left(\alpha_{-h}(a_h) \alpha_{-h}(e_h)\right )\delta_0=a_h\delta_0 \neq 0$, and so we can assume, without loss of generality, that $P_0(r) \neq 0$ (exchanging $r$ for $r \alpha_{-h}(e_h) \delta_{-h}$ if necessary). 

Now, let $J:= \{P_0(s)  :  s \in RrR \ \mbox{and} \ supp(s) \subseteq supp(r) \}$.
Notice that $J$ is a non-empty set that contains $a_0=P_0(r)$ (since $r \in RrR$) and so, since $A$ is $G$-simple, we finish the proof if we show that $J$ is a $G$-invariant ideal of $A$. For this, let $a \in J\cap D_{-h}$. Then $a \delta_0 + \sum_{g} b_g\delta_g \in RrR$ for some  $b_g \in D_g$ and $g \in  supp(r) \setminus \{0\}$. Let $e$ be unit for $a$ in $D_{-h}$. Then $\alpha_h(e) \delta_h  ( a\delta_0 + \sum_g b_g \delta_g)e\delta_{-h} \in RrR$ and
\begin{align*}
\alpha_h(e) \delta_h & ( a\delta_0 + \sum_g b_g \delta_g)e\delta_{-h} \\
 & =\alpha_h(e)\delta_h a \delta_0 e \delta_{-h} + \sum_g \alpha_h(e)\delta_h b_g\delta_g e \delta_{-h}\\
& = \alpha_h(e)\delta_h \alpha_0(\alpha_{-0}(a) e) \delta_{-h} + \sum_g \alpha_h(e)\delta_h(\alpha_g(\alpha_{-g}(b_g)e)\delta_{g-h}\\
& = \alpha_h(e)\delta_h a \delta_{-h} + \sum_g  \alpha_h(\alpha_{-h}(\alpha_h(e))\alpha_g(\alpha_{-g}(b_g)e))\delta_{h+g-h}\\
& = \alpha_h(e a)\delta_0 + \sum_g \alpha_h (e \alpha_g(\alpha_{-g}(b_g)e))\delta_{g}, \ \mbox{since} \ G \ \mbox{is comutative}\\
& = \alpha_h(a)\delta_0 + \sum_g \alpha_h (e \alpha_g(\alpha_{-g}(b_g)e))\delta_{g},  \ \mbox{with} \ g \in supp(r) \setminus \{0\},
\end{align*} 
what implies that $\alpha_h(a) \in J\cap D_h$. Therefore, $J$ is a $G$ - invariant ideal as desired. 

\cqd 

For skew group rings simplicity is related to the center of the ring. Since we are dealing with rings with local units we have to look into corners:

\begin{definition} Let $R=A \rtimes_\alpha G$ be an associative partial skew group ring and $E$ a set of local units for $A$. For each $e\in E$, let $C_e$ be the center of $e\delta_0Re\delta_0$, that is, $C_e := \{x \in e\delta_0Re\delta_0 : xy=yx  \ \forall  y \in e\delta_0Re\delta_0\}.$
\end{definition}

\begin{lema}\label{lema 2}  Suppose we are under the same conditions of lemma \ref{lema 1} and let $e\in E$. Then every non-zero ideal of $R= A \star _{\alpha} G$ has non-empty intersection with $C_e\cap \{e \delta_0 + \sum_{g \in G \setminus \{0\}} b_g\delta_g\}$.
\end{lema}

\demo Let $J$ be a non-zero ideal of $R$ and choose $r \in J \setminus \{0\}$ such that  $\# supp (r)$ is minimal. By lemma \ref{lema 1}, we can find $r'' \in RrR \subseteq J$ such that  $P_0(r'')=e$ and $\# supp (r'') \leq \# supp(r)$. 

Let $r'=e\delta_0r''e\delta_0$. Notice that $r' \in RrR$, $P_0(r')=e$ and  $\# supp (r') \leq \# supp (r'') \leq \# supp(r)$. 

Now, since $P_0(r'')=e$, we have that $P_g(r'e\delta_0a_g\delta_g e\delta_0) = P_g(e\delta_0r''e\delta_0a_g\delta_g e\delta_0) = P_g(e\delta_0a_g\delta_g e\delta_0)= P_g(e\delta_0 a_g \delta_g e\delta_0 r'' e\delta_0)= P_g(e\delta_0a_g\delta_g e\delta_0 r')$. So, since $supp(r'e\delta_0a_g\delta_g e\delta_0)$ and $supp (e \delta_0a_g\delta_g e\delta_0 r')$ are subsets of $\{g\cdot t :  t \in supp (r')\}$ we have that $\# supp(r'e\delta_0a_g\delta_g e\delta_0 - e\delta_0a_g\delta_g e \delta_0r')< \# supp (r') \leq \# supp(r)$, what 
implies that  $r'e\delta_0a_g\delta_g e\delta_0 = e\delta_0 a_g\delta_g e \delta_0r'$ for all $g \in G$ and hence, by linearity, $r' \in C_e$.

\cqd

\begin{theorem}\label{teorema 1} Let $E$ be a set of local units for $A$ and $\alpha=(\{D_t\}_{t\in G}, \{\alpha_t\}_{t \in G} )$ a partial action of an abelian group $G$ such that each ideal $D_t$ has local units. Suppose that $R= A \star _{\alpha} G$ is associative. Then the following are equivalent:
\begin{enumerate}
\item $A \rtimes _{\alpha} G$ is simple.
\item $A$ is $G$-simple and $C_e$ is a field for all $e\in E$. 
\end{enumerate}
\end{theorem}
\demo First suppose that $A \rtimes _{\alpha} G$ is simple. 
The proof that $A$ is $G-$simple is essentially the same as the one done in \cite{Oinert}, one just has to notice that if $J$ is a non-zero proper invariant ideal of $A$ then $(\{J\cap D_t\}_{t\in G}, \{\alpha_t\}_{t \in G} )$   is a partial action and $J \rtimes_\alpha G$ is a non-zero ideal of $A \rtimes_\alpha G$.

Next we show that $C_e$ is a field for all $e\in E$. So, let $a \in C_e$ and consider the ideal generated by $a$ in $A \rtimes_\alpha G$. By the simplicity of  $A \rtimes_\alpha G$, we have that $\langle a \rangle = A \rtimes_\alpha G$.
Thus $e\delta_0 \in \langle a \rangle$ and so there exists $r_i$ e $s_j$ such that
$$
e \delta_0 = \sum_{i,j} r_i a s_j = \sum_{i,j} e\delta_0 r_i a s_j e\delta_0 =  \sum_{i,j} e\delta_0 r_i \underbrace{e\delta_0 a e\delta_0}_{a} s_j e\delta_0, 
$$ and, since $a \in C_e$, we have that
$$
e \delta_0 =  \sum_{i,j}a e\delta_0 r_i e\delta_0  e\delta_0 s_j e\delta_0= a \sum_{i,j} e\delta_0 r_i e\delta_0 s_j e\delta_0. 
$$
So $a$ has an inverse and all we have left to do is show that $a^{-1}\in C_e$. But notice that $a^{-1} = e\delta_0 \left(\sum_{i,j}  r_i e\delta_0 s_j\right) e\delta_0 \in e\delta_0 R e\delta_0$ and, since $a\in C_e$, we have that for all $y\in e\delta_0 R e\delta_0$ it holds that 
$ay=ya \Rightarrow  e\delta_0y=a^{-1}ya \Rightarrow  ya^{-1} = a^{-1}yaa^{-1} 
\Rightarrow  ya^{-1} = a^{-1}y$. We conclude that $a^{-1} \in C_e$ and hence $C_e$ is field.

Now suppose that $A$ is $G$-simple and $C_e$ is a field for all $e\in E$. Let $J$ be a nonzero ideal of  $R=A \rtimes_\alpha G$. By Lemma \ref{lema 2} there is a non-zero $r \in J\cap C_e$ for every $e\in E$. Since $C_e$ is field this implies that  $e\delta_0= r^{-1}r \in J $ for all $e\in E$ and, since $E\delta_0$ is a set of local units for $R$, we conclude that $J=R$.
\cqd

\section{An Application to set dynamics}

In \cite{BG} it was shown that there is a one to one correspondence between partial actions in a set $X$ and partial actions in $\Fo$, where $\Fo$ is the algebra of all functions from $X$ to a field $K$ with finite support (see \cite{BG}). More precisely, if $\theta=(\{X_t\}_{t \in G}, \{h_t\}_{t \in G})$ is a partial action in $X$ then $\alpha=(\{D_t\}_{t \in G}, \{\alpha_t\}_{t \in G})$, where $D_t=\mathcal{F}_0(X_t)=\{f\in \Fo: \ f(x)=0 \ \forall \ x\notin X_t \}$ and $\alpha_t(f):=f\circ h_{-t}$,
is a partial action of $G$ in $\Fo$.

Our goal in this section is to show the following theorem:
\begin{theorem}\label{teorema2} Let $G$ be an abelian group and $\theta=(\{X_t\}_{t \in G}, \{h_t\}_{t \in G})$ a partial action in a set $X$. Then $\Fo \star G$ is simple if, and only if, $\theta$ is a minimal and free partial action.
\end{theorem}

Of course we will use theorem \ref{teorema 1} to prove the above result. So we have to check that the hypotheses are verified. Notice that $\Fo \star G$ is associative (see \cite{BG}) and it is clear that $\Fo$, as well as the ideals $F(X_t)$, have local units and so we can apply theorem \ref{teorema 1}. But theorem \ref{teorema 1} also requires that  we choose a set of local units for $\Fo$. So we let $E:=\{\chi_A: \ A \text{ is a finite subset of } X\}$, where $\chi_A$ denotes the characteristic function of $A$, be a fixed set of local units for $\Fo$.


We now recalll the relevant definitions mentioned in theorem \ref{teorema2}.


\begin{definition}
A partial action $\theta=(\{X_t\}_{t \in G}, \{h_t\}_{t \in G})$ of a group $G$ in $X$ is minimal if the only $G$-invariant subsets of $X$ are $\emptyset$ and $X$. $\theta$ is free if, for all $x\in X$, $h_t(x)= x$ implies that $t=0$.
\end{definition}

\begin{obs} Minimality can also be characterized in other ways, more specifically, $\theta$ is minimal iff for all $x \in X$, $V_x =\{ h_t(x)  : \ t \in G \ \mbox{and} \ x \in X_{-t} \} = X$, what is equivalent to say that if $U$ and $V$ are subsets of $X$, then there exists $t \in G$ such that $h_t(U)\cap V \neq \emptyset$.
\end{obs}

Part of theorem \ref{teorema2} follows from the correspondence between $G$-invariant sets of $X$ and $G$-invariant sets of $\Fo$:

\begin{proposition}
Let $\theta=(\{X_t\}_{t \in G}, \{h_t\}_{t \in G})$ be a partial action in a set $X$.  Then, 
$V \subset X$ is $G$-invariant if and only if $\mathcal{F}_0(V)$ is $G$-invariant.
\end{proposition}
\demo 
First suppose that $V$ is $G$-invariant and let $f \in \mathcal{F}_0(V)\cap D_{-t}$. Then $\alpha_t(f)(x)\neq 0$ implies that $h_{-t}(x)\in V \cap X_{-t} $ and hence $x \in h_t(V\cap X_{-t})\subseteq V$. So $\alpha_t(f) \in \mathcal{F}_0(V)$.

Now suppose that $\mathcal{F}_0(V)$ is $G$-invariant and assume that there exists a $x \in V \cap X_{-t}$ such that $h_t(x) \notin V $. Notice that $\delta_x \in \mathcal{F}_0(V) \cap \mathcal{F}_0(X_{-t})$ and hence $\alpha_t(\delta_x)\in \mathcal{F}_0(V) $. But $\alpha_t(\delta_x)(h_t(x))=\delta_x\circ h_{-t}(h_t(x)) = \delta_x(x)=1$, a contradiction. So $h_t(x) \in V$ and $V$ is $G$-invariant.
\cqd

\begin{corolario} A partial action $\theta=(\{X_t\}_{t \in G}, \{h_t\}_{t \in G})$ in a set $X$ is minimal if, and only if, $\Fo$ is $G$-simple.
\end{corolario}

Next we will show that under the additional hypothesis that $\theta$ is free then, for each $e\in E$, $C_e$ is a field.

\begin{lema}\label{CecontidoF0}
Let $\theta$ be a free partial action of an abelian group $G$. Then for all $e\in E$, say $e=\chi_A$, $C_e \subset e \Fo e = \mathcal{F}_0(A)$. 
\end{lema}
\demo 
Suppose that $x=e\delta_0 \left( \sum_g f_g \delta_g \right)e\delta_0 \in C_e$ and there exists $g\neq 0$ such that $z:=e\delta_0 f_g \delta_g e \delta_0 \neq 0$. We will derive a contradiction.

Notice that, for all $t\in G$ and $f_t \in D_t$, since $G$ is abelian and $x\in C_e$, we have that 
$z( e \delta_0 f_t \delta_t e \delta_0)  = P_{g+t}\left(x(e \delta_0 f_t \delta_t e \delta_0)\right)\delta_{g+t} = P_{g+t}\left((e \delta_0 f_t \delta_t e \delta_0)x\right)\delta_{g+t}= ( e \delta_0 f_t \delta_t e \delta_0)z,$
and hence $ z \in C_e$. So $\alpha_g(\alpha_{-g}(ef_g)ef_0)\delta_g = z e\delta_0 f_0 \delta_0 e \delta_0 = 
 e\delta_0 f_0 \delta_0 e \delta_0 z =  \alpha_g(\alpha_{-g}(f_0 e f_g)e)\delta_g,$ for all $f_0\in \Fo$ and, since $\alpha_g$ is an isomorphism, this implies that $\alpha_{-g}(ef_g)ef_0 = \alpha_{-g}(f_0 e f_g) e $ for all $f_0\in \Fo$, which is equivalent to 
\begin{equation}\label{eq1}
 f_g|_A (h_g(x)). \chi_A(x).f_0(x)  = f_0(h_g(x)).f_g |_A(h_g(x)). \chi_A(x) 
\end{equation}
 for all $f_0\in \Fo$ and $x\in X_{-g}$.

Now, $z \neq 0$ implies that $e f_g \neq 0$ and so $f_g|_A \neq 0$. Furthermore, $e f_g \delta_g e \delta_0 = \alpha_g(\alpha_{-g}(ef_g) e )\delta_g \ne 0 $ and so $\alpha_{-g}(ef_g) e \neq 0$. Therefore there exists $x\in X_{-g}$ such that $f_g|_A(h_g(x)).\chi_A(x) \neq 0$. Let $f_0 = \chi_{\{ x \} }$. Then the left side of equation \ref{eq1} is nonzero and, since the action is free, the right side is zero, a contradiction. 
\cqd

\begin{proposition}\label{CeiguF0}
If $\theta = (\{X_t\}_{t \in G}, \{h_t\}_{t \in G})$ is a free, minimal partial action of an abelian group in a set $X$ then, for all $e\in E$, $C_e$ is a field. More precisely, if $e= \chi_A\in E $, then $C_e\subseteq \{f \in \Fo \ :  \ supp(f)  =A \} \cup \{0\}$. 
\end{proposition}

\demo
Let $f \in C_e\subseteq \mathcal{F}_0(A)$ be a non-zero function. Then, for all $g\in G$ and $f_g \in D_g$, we have that
$f(e\delta_0f_g\delta_ge\delta_0) =(e\delta_0f_g\delta_ge\delta_0)f$, so $\alpha_g(\alpha_{-g}(f . f_g)\chi_A)\delta_g = \alpha_g(\alpha_{-g}(\chi_A . f_g)f)\delta_g,$  for all $g\in G$ and $f_g \in D_g$ and hence
\begin{equation}\label{equation2}
f(x)f_g(x)\chi_A(h_{-g}(x)) = \chi_A(x)f_g(x)f(h_{-g}(x))\ \ \forall g\in G, \ f_g \in D_g \text{ and } x \in X_g.
\end{equation}
Now suppose that $supp(f) \subsetneq A$ and let $y \in supp(f).$ Since $\theta$ is minimal there exists $t \in G$ such that $h_{-t}(y) \in A \setminus supp(f)$. So equation \ref{equation2}, with $g=t$ and $f_g = \delta_y$, becomes $f(x)\delta_y(x) \chi_A(h_{-t}(x))=\chi_A(x) \delta_y(x)f(h_{-t}(x)) \ \forall x\in X_t$, and hence, for $x=y$, we have that $f(y)=0$, a contradiction. We conclude that $supp(f)=A$ and so there exists $f^{-1}$ such that $f f^{-1} = f^{-1} f = e$. The proof that $f^{-1} \in C_e$ is analogous to what was done in the proof of theorem \ref{teorema 1}. 
\cqd

The following proposition proves the last part of theorem \ref{teorema2}.

\begin{proposition}\label{simple-livre}
 If $\Fo \rtimes G$ is simple then $\theta=(\{X_t\}_{t \in G}, \{h_t\}_{t \in G})$ is free.
 \end{proposition}
\demo
Suppose that $\theta$ is not free. Then there exists $x\in X$ and $g \in G$, $g \neq 0$, such that $x \in X_{-g}$ and $h_g(x)=x$. Consider the ideal $I$ generated by $\chi_x\delta_0 - \chi_x \delta_g$ (notice that $\chi_x \in \mathcal{F}_0(X_g)$ since $x=h_g(x)\in X_g$). 

We will show that the sum of coefficients of elements in $I$ is zero. For this, notice that $\alpha_{-g}(\chi_x) = \chi_x$ and so $$
a_s\delta_s (\chi_x\delta_0 -\chi_x \delta g) b_t \delta_t = a_s \delta_s \chi_x b_t \delta_t - a_s \delta_s \alpha_g( \alpha_{-g}(\chi_x) b_t)\delta_{t+g} =
a_s \delta_s \chi_x b_t \delta_t - a_s \delta_s \alpha_g(\chi_x b_t) \delta_{t+g}.  $$ 

Now, $\alpha_g(\chi_x b_t) \neq 0$ if, and only if, there exists $y\in X$ such that $\chi_x ( h_{-g}(y)) b_t (h_{-g}(y)) \neq 0$ and this is true if, and only if,  $b_t (h_{-g}(y)) \neq 0$ and $h_{-g}(y)=x$, that is, $y=h_g(x)=x$, in which case $\alpha_g(\chi_x b_t)(x) =  b_t(x)$. So $\alpha_g(\chi_x b_t) = \chi_x b_t$ and hence
$$
a_s\delta_s (\chi_x \delta_0 -\chi_x \delta g) b_t \delta_t  =  a_s \delta_s \chi_x b_t \delta_t - a_s \delta_s \chi_x b_t \delta_{t+g} =\alpha_s(\alpha_{-s}(a_s)\chi_x b_t) \delta_{t+s} - \alpha_s(\alpha_{-s}(a_s)\chi_x b_t) \delta_{t+g+s}
.$$

We conclude that the sum of coefficients of elements in $I$ is zero. But then $\chi_x \delta_0 \notin I$ and hence $\Fo \rtimes G $ is not simple.
\cqd

We finalize this section with an example of a minimal, free partial action of the group of the integer numbers, denoted by $\Z$, in the set of natural numbers, denoted by $\N$.

\begin{example}  Let $X_0 = \N$, $h_0=id$, $X_{-1} = \N$, $X_1= \N - \{1\}$ and $h:X_{-1}\rightarrow X_1$ be defined by $h(n) = n+1$. For all other $t \in \Z$ let $X_{-t}$ be the domain of $h^t$ and $h_t = h^t$. Then $\{(X_t,h_t)\}$ is a free, minimal partial action and hence the associated partial skew group ring $\Fo\star \Z$ is simple.
\end{example}

\section{An application to topological dynamics}

We now turn our attention to the context of topological partial actions. In this setting the correspondence between partial actions in a locally compact Hausdorff space $X$ and partial actions in the C*-algebra of continuous functions vanishing at infinity, $\mathcal{C}_0(X)$, is well known (see \cite{ELQ} for example) and follows the same ideas we present in the previous section, namely, if $\theta=(\{X_t\}_{t \in G}, \{h_t\}_{t \in G})$ is a partial action in $X$ then $\alpha=(\{D_t\}_{t \in G}, \{\alpha_t\}_{t \in G})$, where $D_t=\mathcal{C}_0(X_t)$ and $\alpha_t(f):=f\circ h_{-t}$, is a partial action of $G$ in $\mathcal{C}_0(X)$. Simplicity of the associated C*-partial crossed product is studied in \cite{ELQ}, where it is shown that if the action is topologically free and minimal then the associated partial crossed product is simple. Minimality of a topological action is exactly what one expects, that is, there are no proper, open invariant subsets, what is equivalent to say that the orbits are dense. We recall the definition of topological freeness below.

\begin{definition} A topological partial action $\theta=(\{X_t\}_{t \in G}, \{h_t\}_{t \in G})$ is topologically free if for all $t \neq 0$ the set $F_t=\{ x \in X_{-t}: h_t(x) = x \}$ has empty interior.
\end{definition}

Using theorem \ref{teorema 1} we will show the following:

\begin{theorem}\label{teorema3} Let $\theta=(\{X_t\}_{t \in G}, \{h_t\}_{t \in G})$ be a partial action of an abelian group in a compact space $X$ such that each $X_t$ is a clopen set. Then the partial skew group ring $\mathcal{C}(X) \star G $ is simple if, and only if, $\theta$ is topologically free and minimal. 
\end{theorem}

\begin{obs} Partial actions on the Cantor set by clopen subsets are exactly the ones for which the enveloping space is Hausdorff, see \cite{EGG} 
\end{obs}  

The proof of the above theorem will follow the same ideas presented in the previous section. Actually the proofs just need to be adapted to the case in hand. We show the relevant details below, but before we proceed notice that we can apply theorem \ref{teorema 1} to prove the result above, since by \cite{Ex} the partial skew group ring is associative, and since the partial action acts on clopen sets, each $D_t$ has a unit. Furthermore the ring $\mathcal{C}(X)$ is unital and hence the set of local units required in theorem \ref{teorema 1} may be taken as the unit in $\mathcal{C}(X)$, which we denote by 1. 

\begin{proposition} A partial action $\theta=(\{X_t\}_{t \in G}, \{h_t\}_{t \in G})$ in a compact space $X$ is minimal if, and only if, $\mathcal{C}(X)$ is $G$-simple.
\end{proposition}
\demo
The proof of this can be found in \cite{ELQ}.
\cqd

\begin{lema}
Let $\theta$ be a topologically free partial action of an abelian group $G$. Then $C_1 \subset  \mathcal{C}(X)$. 
\end{lema}
\demo 
Suppose that $x=\sum_t f_t \delta_t  \in C_1$ and there exists $g\neq 0$ such that $ f_g \neq 0$. Notice that the first part of the proof of lemma \ref{CecontidoF0} was done in general, so that in the case at hand equation \ref{eq1} reduces to 
\begin{equation}\label{eq1a}
f_g(h_g(x)) f_0(x) = f_g(h_g(x)) f_0(h_g(x)),
\end{equation}
for all $f_0 \in \mathcal{C}(X)$ and $x \in X_{-g}$.

Now, since $f_g \neq 0$, we have that $\alpha_{-g}(f_g) \neq 0$ and so there exists $x \in X_{-g}$ such that $f_g(h_g(x)) \neq 0$. Since $f_g$ is continuous there exists an open neighborhood $V$ of $x$ such that $f_g(h_g(y)) \neq 0$ for all $y\in V$. Consider the open neighborhood  $V\cap X_{-g}$ of $x$. Since $\theta$ is topologically free there exists $y \in V\cap X_{-g}$ such that $h_g(y)\neq y$. Then, by Urysohn´s lemma, there exists $f_0 \in \mathcal{C}(X)$ such that $f(y)=1$ and $f(h_g(y))=0$. But then, for this $f_0$ and $y$ equation \ref{eq1a} leads to a contradiction.
\cqd

\begin{proposition}
If $\theta = (\{X_t\}_{t \in G}, \{h_t\}_{t \in G})$ is a topologically free, minimal partial action of an abelian group in a compact space $X$ then $C_1$ is a field. More precisely, $C_1 = \C \cdot 1$, that is, $C_1$ is the algebra of constant functions. 
\end{proposition}

\demo
Let $f \in C_1\subseteq \mathcal{C}_0(X)$ be a non-zero function. Notice that the first part of proposition \ref{CeiguF0} was done in general and so it is valid in the case at hand, for which equation \ref{equation2} becomes 
$$ f(x)f_g(x) = f_g(x)f(h_{-g}(x))\ \ \forall g\in G, \ f_g \in D_g \text{ and } x \in X_g.$$

Now for each $g\in G$ let $f_g$ be the unit for $D_g$, that is, $f_g = \chi_{X_g}$. Then the above equation implies that $f(x) = f(h_{-g}(x))$ for all $x \in X_g$ and $g\in G$ and, since $\theta$ is minimal and $f$ is continuous, we obtain that $f$ is constant as desired.
\cqd

The following will finish the proof of theorem \ref{teorema3}.

\begin{proposition}
 If $\mathrm{C}(X) \rtimes G$ is simple then $\theta=(\{X_t\}_{t \in G}, \{h_t\}_{t \in G})$ is topologically free.
 \end{proposition}
\demo
Suppose that $\theta$ is not topologically free. Then there exists $g\neq 0$ in $G$ such that the interior of $F_g$ is not empty. Let $x$ be an element in the interior of $F_g$. By Urysohn´s lemma there exists a continuous function $f$ such that $f(x) = 1$ and the support of $f$ is contained in the interior of $F_g$. 

Notice that $f = \chi_{F_g} \cdot f$ and hence $\alpha_g(f ) =f = \alpha_{-g}(f)$. Now consider the ideal generated by $f \delta_0 - f \delta_g$. Proceeding analogously to what was done in proposition \ref{simple-livre}, that is, expanding terms of the form $a_s \delta_s(f \delta_0 - f\delta_g) b_t \delta_t$, we have that the sum of coefficients of elements in $I$ is zero. 
But then $f \delta_0 \notin I$ and hence $\Fo \rtimes G $ is not simple.
\cqd

\vspace{1.5pc}

{\bf Acknowledgement:}: I would like to thank Viviane Beuter for the useful discussions on the topic. Work was partially supported by CNPq.

\vspace{0.5pc}

D. Goncalves, Departamento de Matemática, Universidade Federal de Santa Catarina, Florianópolis, 88040-900, Brasil

Email: daemig@gmail.com

\vspace{0.5pc}

\end{document}